\documentclass[a4paper]{amsart}
\usepackage{amssymb}
\usepackage{mathtools}
\usepackage[alphabetic,nobysame,non-compressed-cites,initials]{amsrefs}
\usepackage[hidelinks]{hyperref}
\usepackage{graphicx}
\usepackage{todonotes}
\usepackage{enumitem}
\usepackage{CJKutf8}

\newcommand{\R}{\mathbb{R}}
\newcommand{\C}{\mathbb{C}}
\newcommand{\T}{\mathbb{T}}

\newcommand{\tphi}{\tilde\phi}
\newcommand{\tQ}{\tilde Q}
\newcommand{\ii}{{\rm i}}

\newcommand{\re}{\operatorname{Re}}
\newcommand{\im}{\operatorname{Im}}
\newcommand{\res}{\operatorname{Res}}

\theoremstyle{plain} 
\newtheorem{theorem}{Theorem}

\newtheorem{proposition}[theorem]{Proposition}
\newtheorem*{problem*}{Problem}
\theoremstyle{definition} 

\theoremstyle{remark}
\newtheorem{remark}{Remark}
\newtheorem{example}{Example}

\begin{document}

\title{Helicoids and vortices}

\author{Hao Chen}
\address[Chen]{Georg-August-Universit\"at G\"ottingen, Institut f\"ur Numerische und Angewandte Mathematik, 37083 G\"ottingen, Germany}
\address[Chen]{ShanghaiTech University, Institute of Mathematical Sciences, 201210 Shanghai, P.R.\ China}
\email{hao.chen.math@gmail.com}
\thanks{H.\ Chen is partially supported by Individual Research Grant from Deutsche Forschungsgemeinschaft within the project ``Defects in Triply Periodic Minimal Surfaces'', Projektnummer 398759432.}

\author{Daniel Freese}
\address[Freese]{Department of Mathematics, Indiana University, Bloomington, IN 47405, USA}
\email{djfreese@iu.edu}

\keywords{minimal surfaces, vortex dynamics}
\subjclass[2020]{Primary 53A10, 76B47}

\date{\today}

\begin{abstract}
	We point out an interesting connection between fluid dynamics and minimal
	surface theory: When gluing helicoids into a minimal surface, the limit
	positions of the helicoids correspond to a ``vortex crystal'', an equilibrium
	of point vortices in 2D fluid that move together as a rigid body.  While
	vortex crystals have been studied for almost 150 years, the gluing
	construction of minimal surfaces is relatively new.  As a consequence of the
	connection, we obtain many new minimal surfaces and some new vortex crystals
	by simply comparing notes.
\end{abstract}

\maketitle

In 2005,
Traizet and Weber~\cite{traizet2005} glued helicoids into screw-motion
invariant minimal surfaces with helicoidal ends.  For the glue construction to
succeed, the limit positions of the helicoids must satisfy a balancing
condition and a nondegenerate condition.  For simplicity, they assumed that the
helicoids are aligned along a straight line, and noticed that the roots of
Hermite polynomials provide examples of balanced and nondegenerate
configurations.  Recently, the second named author implemented a similar
construction without this assumption~\cite{freese2021}.

The main goal of this note is to point out an interdisciplinary connection: The
balanced configurations of helicoids correspond to binary vortex crystals.
Here, a \emph{vortex crystal}~\cite{aref2003}, also known as \emph{vortex
equilibrium}, is a configuration of vortices in 2D fluids that moves as a rigid
body, i.e.\ without any change of shape and size.  A vortex crystal is
\emph{binary} if the circulations of vortices are $\pm 1$.  We will recall
vortex dynamics in Section~\ref{sec:vortex}.  For an example of this
connection: the ``definite'' configurations in~\cite{traizet2005}, given by
roots of Hermite polynomials, correspond to vortex crystals that trace back to
19th century~\cite{stieltjes1885}. 

Rotating vortex crystals correspond to screw-motion invariant minimal surfaces;
this connection was readily established in~\cites{traizet2005, freese2021}.
Our main results in Section~\ref{sec:mainresult} establish the other cases of
the claimed connection.  The construction will be given in
Section~\ref{sec:construct}, where we glue helicoids into translation-invariant
minimal surfaces, corresponding to translating or periodic stationary vortex
crystals.  The construction will only be sketched because similar constructions
have been repeated many times in the literature~\cites{traizet2008,
traizet2008b, chen2021}.

In view of the connection, we will compare notes between fluid dynamics and
minimal surface theory, and obtain new examples in Section~\ref{sec:example}
for both vortex crystals and minimal surfaces.

The minimal surface theory would benefit a lot because, in about 150 years, the
fluid dynamics community has accumulated a large collection of examples of
binary vortex crystals.  In particular, stationary and translating vortex
crystals have been analytically obtained with the help of Adler--Moser
polynomials.  When symmetries are imposed, nondegeneracy was recently verified
for translating Adler--Moser examples~\cite{liu2020}.  They then lead to many
examples.

On the other hand, the idea to glue helicoids into minimal surfaces is only
about 18 years old.  Nevertheless the minimal surface community is also in
possession of a few examples that would lead to new vortex crystals.  For
instance, it is very hard to observe periodic vortex crystals in experiment,
hence there has been very few relevant physics studies; but our knowledge on
periodic minimal surfaces is quite advanced~\cites{meeks1988, meeks1990}.

\medskip

It is worth mentioning that, in natural sciences, minimal surfaces have been
used to model bicontinuous structures in liquid crystals, block copolymers,
biological membranes, and nuclear pasta, etc.; see e.g.~\cite{hyde1996}.  In
particular, helicoid structures are used to model topological defects (screw
dislocations) in layered structures~\cites{kamien2006, santangelo2007,
matsumoto2012, terasaki2013, horowitz2015, schneider2016, dasilva2021}.  On the
other hand, analogy between these structures and vortices in superfluids and
superconductors has been noticed since long~\cites{degennes1972, lubensky1990}.

Nevertheless, rigorous mathematical treatment of helicoidal minimal surfaces
has been missing until Traizet developed the node-opening
technique~\cites{traizet2002, traizet2002b, traizet2008, traizet2008b,
traizet2005} for constructing minimal surfaces.  The present work is the first
to state a well-defined mathematical analogy between point vortices and
helicoids.  We are able to push the analogy to its full generality and to
provide rigorous constructions using the node-opening technique.

\begin{remark}
	Our work is not the first to connect minimal surfaces with fluid dynamics.
	In~\cite{traizet2015}, a correspondence was established between hollow
	vortices and minimal surfaces bounded by horizontal symmetry curves.
\end{remark}

\subsection*{Acknowledgement}

H.\ Chen is grateful to Prof.\ Yong Liu for inspiring discussions.

\section{Vortex dynamics}\label{sec:vortex}

We recommend~\cite{aref2003} for general reference on vortex crystals.

Incompressible and inviscid flow in zero gravity is governed by the Euler
equation
\begin{equation}\label{eq:euleru}
	\big( \frac{\partial}{\partial t} + u \cdot \nabla \big) u = - \nabla p
\end{equation}
under the incompressible condition $\nabla \cdot u = 0$, where $u$ is the flow
velocity field and $p$ is the pressure.  The pressure term can be eliminated by
taking the curl of the Euler equation, which results in the reformulation
\begin{equation}\label{eq:eulerw}
	\big( \frac{\partial}{\partial t} + u \cdot \nabla \big) w
	= w \cdot \nabla u
\end{equation}
in terms of the vorticity field $w := \nabla \times u$.  For 2-dimensional
flows, the right-hand side of~\eqref{eq:eulerw} vanishes, giving the transport
equation
\[
	\big( \frac{\partial}{\partial t} + u \cdot \nabla \big) w = 0.
\]
We see that vorticity is transported in the velocity field as material
elements.

A point vortex in the plane is given by the vorticity field
\[
	w(x) = \frac{\sigma}{2\pi}\delta(x),
\]
where $\delta$ is the Dirac delta function, and $\sigma \in \R$ is the
\emph{circulation} of the vortex.  In the complex coordinate, a point vortex
gives rise to a velocity field
\[
	\overline{u(z)} = \frac{1}{2\pi\ii} \frac{\sigma}{z}.
\]
We use $(p_k, \sigma_k)_{1 \le k \le n}$ to denote a configuration of $n$ point
vortices located at $p_k$ with circulation $\sigma_k$, $k = 1, \cdots, n$.
Each vortex is advected as a material particle by the velocity field produced
by other vortices.  So the dynamics of the vortex configuration is governed by
the ordinary differential equation
\begin{equation}\label{eq:vortexdynamics}
	\frac{d}{dt} \overline{p_j} = \frac{1}{2\pi\ii} \sum_{k \ne j}
	\frac{\sigma_k}{p_j-p_k},
	\quad \forall 1 \le j \le n,
\end{equation}
also known as the Kirchhoff--Routh dynamics~\cites{lin1941, davila2020}.

We say that the configuration $(p_k, \sigma_k)_{1 \le k \le n}$ is a
\emph{vortex crystal} if it moves as a rigid body.  If this is the case, we
have $d p_j / dt = v + \ii \omega p_j$, where $v \in \C$ and $\omega \in \R$
are constant for all vortices.  Hence vortex crystals are characterized by the
balance equations
\begin{equation}\label{eq:balancevortex}
	F_j := - \overline{v} + \ii \omega \overline{p_j} +
	\frac{1}{2\pi\ii} \sum_{k \ne j} \frac{\sigma_k}{p_j-p_k}=0,
	\quad \forall 1 \le j \le n.
\end{equation}
Multiply~\eqref{eq:balancevortex} by $\sigma_j$, and sum the conjugates over
$j$, we obtain
\begin{equation}\label{eq:moment1}
	v \sum_{j=1}^n \sigma_j + \ii \omega \sum_{j=1}^n \sigma_j p_j = 0.
\end{equation}
Multiply~\eqref{eq:balancevortex} by $\sigma_j p_j$, and take the sum over $j$,
we obtain
\begin{equation}\label{eq:moment2}
	\overline{v} \sum_{j=1}^n \sigma_j p_j
	- \ii \omega \sum_{j=1}^n \sigma_j |p_j|^2 =
	\frac{1}{4\pi\ii} \bigg[
		\Big( \sum_{k=1}^n \sigma_k \Big)^2 - \sum_{k=1}^n \sigma_k^2
	\bigg].
\end{equation}

In this paper, we only consider binary vortex crystals, so $\sigma_k = \pm 1$.
Let $n_\pm$ be, respectively, the number of vortices with circulation $\pm 1$,
and write $m = n_+ - n_- = \sum \sigma_k$.  We distinguish three cases,
\begin{enumerate}
	\item  We say that the vortex crystal is \emph{rotating} if $\omega \ne 0$.
		In this case, the governing equation~\eqref{eq:balancevortex} is invariant
		under Euclidean rotations.  Moreover, we may assume that $v=0$ up to a
		translation, and that $\omega=1$ up to a Euclidean scaling.  After these
		normalizations, \eqref{eq:moment1} and~\eqref{eq:moment2} give
		\[
			\sum_{k=1}^n \sigma_k p_k = 0 \quad \text{and} \quad
			\sum_{k = 1}^n \sigma_k |p_k|^2 = \frac{m^2 - n}{4\pi}.
		\]

	\item We say that the vortex crystal is \emph{translating} if $\omega=0$ but
		$v \ne 0$.  Then~\eqref{eq:moment1} implies that $m=0$ so $n_+ = n_- =
		n/2$.  In this case, \eqref{eq:balancevortex} is
		invariant under translations.  We may assume that $v = 1$ up to a complex
		scaling (Euclidean scaling and rotation).  After this normalization,
		\eqref{eq:moment2} implies that
		\[
			\sum_{k=1}^n \sigma_k p_k = -\frac{n}{4\pi\ii}.
		\]

	\item We say that the vortex crystal is \emph{stationary} if $\omega=0$ and
		$v=0$.  Then~\eqref{eq:moment2} implies that $m^2 = n$, hence $n_+$ and
		$n_-$ must be successive triangular numbers.  In this case,
		\eqref{eq:balancevortex} is invariant under Euclidean similarities
		(translations, rotations, and scalings).
\end{enumerate}

A vortex crystal is said to be \emph{stable} if any sufficiently small
perturbation does not diverge.  Linear and nonlinear analyses have been carried
out on the stability.  We say that a vortex crystal is
\emph{nondegenerate\footnotemark} if the Jacobian matrix $DF = (\frac{\partial
F_i}{\partial p_j})_{i,j}$ for the autonomous system~\eqref{eq:vortexdynamics}
has the maximum possible rank.  We have seen that, for a rotating (resp.\
translating, stationary) vortex crystal, the dynamics is invariant under
Euclidean rotations (resp.\ translations, similarities), so the maximum
possible real rank of its Jacobian is $2n-1$ (resp. $2n-2$, $2n-4$).

\footnotetext{In fluid dynamics literature, degenerate vortex crystals were
often said to be ``neutrally stable''.}

\medskip

We may also consider singly or doubly periodic vortex crystals.  In this case,
it is convenient to consider a co-rotating reference frame in which the periods
are fixed.  Then the vortex crystal is either translating or stationary.

\begin{itemize}
	\item Assume that the vortex crystal is \emph{singly periodic}, i.e.\
		invariant under a single translation $T \in \C$.  Up to rotations and
		scalings, we may fix $T = 1$.  Then the vortices can be seen as lying in
		the annulus $\C/\langle 1 \rangle$.  Up to translations, a vortex crystal
		in the annulus is governed by
		\[
			F_j := - \overline{v} +
			\frac{1}{2\pi\ii} \sum_{k \ne j} \pi\sigma_k\cot\pi(p_j-p_k) = 0,
			\quad \forall 1 \le j \le n
		\]
		in an appropriate reference frame.

	\item Assume that the vortex crystal is \emph{doubly periodic}, i.e.\
		invariant under two linearly independent translations $T_1$ and $T_2$.  Up
		to rotations and scalings, we may fix $T_1 = 1$ and $T_2 = \tau \in \C$.
		Then the vortices can be seen as lying in the flat torus $\C/\langle 1,
		\tau \rangle$.  Up to translations, a vortex crystal in the torus is
		governed by
		\[
			F_j := - \overline{v} + \frac{1}{2\pi\ii} \sum_{k \ne j} \sigma_k
			\big(\zeta(p_j-p_k; \tau) - \xi(p_j-p_k; \tau)\big) = 0,
			\quad \forall 1 \le j \le n
		\]
		in an appropriate reference frame, where $\zeta(z;\tau)$ is the Weierstrass
		zeta function on the torus $\C/\langle 1, \tau \rangle$ and $\xi(z; \tau) =
		2 x \zeta(1/2; \tau) + 2 y \zeta(\tau/2; \tau)$ with $z = x + y \tau$, $x,
		y \in \R$.
\end{itemize}

In either case, the maximum possible real rank of the Jacobian is $2n-2$, which
defines the degeneracy of these vortex crystals.

\section{Main results}\label{sec:mainresult}

The connection between rotating vortex crystals and screw-motion invariant
minimal surfaces was already established in the following
theorem~\cite{freese2021}.

\begin{theorem}[Rotating vortex crystal]\label{thm:screw}
	Let $(p^\circ_k, \sigma_k)_{1 \le k \le n}$ be a normalized nondegenerate
	rotating vortex crystal with $n$ vortices at $p^\circ_1$, \ldots, $p^\circ_n
	\in \C$.  Then there exists a one-parameter family $(M_\varepsilon)_{0 <
	\varepsilon < \delta}$ of embedded minimal surfaces in $\R^3$ such that:
	\begin{enumerate}
		\item $M_\varepsilon$ admits a screw symmetry $S_\varepsilon$ composed of a
			vertical translation $2\pi(0, 0, 1)$ and a rotation around the vertical
			axis by an angle $2\pi\varepsilon^2$.

		\item The quotient $M_\varepsilon / S_\varepsilon$ is of genus $n-1$ and
			has two ends.  The ends are helicoidal if $m = \sum\sigma_k \ne 0$, or
			planar if $m = 0$.

		\item As $\varepsilon \to 0$, up to a translation, $M_\varepsilon$
			converges to a helicoid of period $(0,0,2\pi)$ in the neighborhood of
			$(p^\circ_i/\varepsilon, 0)$ for each $1 \le i \le n$.  The helicoid is
			right-handed (resp.\ left-handed) if $\sigma_i = 1$ (resp.\ $-1$).  

		\item After rescaling the horizontal coordinates by $\varepsilon$, the
			resulting surface (no longer minimal) converges to the union of the
			multigraph of the multivalued function
			\[
				f(z) = \sum_{i = 1}^n \sigma_i \arg(z-p^\circ_i), \qquad z \in \C - \{p^\circ_1, \cdots, p^\circ_n\},
			\]
			the multigraph of $f(z)+\pi$, and vertical lines over the points $p^\circ_i$.
	\end{enumerate}
	Conversely, for any such family of screw-motion invariant minimal surfaces,
	in the limit $\varepsilon\to 0$, the positions and handednesses of the
	helicoids give a rotating vortex crystal (possibly degenerate).
\end{theorem}

So we only need to construct translation invariant minimal surfaces,
corresponding to translating or stationary vortex crystals.  More specifically,
we will prove the theorems below.

\begin{theorem}[Finite translating vortex crystals]\label{thm:translate}
	Let $(p^\circ_k, \sigma_k)_{1 \le k \le n}$ be a nondegenerate translating
	vortex crystal normalized with velocity $v = 1$, with $n$ vortices at
	$p^\circ_1$, \ldots, $p^\circ_n \in \C$.  Then there exists a one-parameter
	family $(M_\varepsilon)_{0 < \varepsilon < \delta}$ of embedded minimal
	surfaces in $\R^3$ such that:
	\begin{enumerate}
		\item $M_\varepsilon$ admits a translational symmetry $T_{0,\varepsilon} =
			2\pi(-2\pi\varepsilon, 0, 1)$.  So $M_\varepsilon$ is a singly periodic
			minimal surface in $\R^3$.

		\item The quotient $M_\varepsilon / T_{0,\varepsilon}$ is of genus $n-1$
			and has two helicoidal ends if the vortex crystal is translating;

		\item As $\varepsilon \to 0$, up to a translation, $M_\varepsilon$
			converges to a helicoid of period $(0,0,2\pi)$ in the neighborhood of
			$(p^\circ_i/\varepsilon, 0)$ for each $1 \le i \le n$.  The helicoid is
			right-handed (resp.\ left-handed) if $\sigma_i = 1$ (resp.\ $-1$).  

		\item After rescaling the horizontal coordinates by $\varepsilon$, the
			resulting surface (no longer minimal) converges to the union of the
			multigraph of the multivalued function
			\[
				f(z) = \sum_{i = 1}^n \sigma_i \arg(z-p^\circ_i), \qquad z \in \C - \{p^\circ_1, \cdots, p^\circ_n\},
			\]
			the multigraph of $f(z)+\pi$, and vertical lines over the points $p^\circ_i$.
	\end{enumerate}
	Conversely, for any such family of singly periodic minimal surfaces, in the
	limit $\varepsilon\to 0$, the positions and handednesses of the helicoids
	give a translating vortex crystal (possibly degenerate).
\end{theorem}

\begin{remark}
	The helicoid limits of minimal surfaces are necessarily balanced, hence
	correspond to vortex crystals.  In the opposite direction, for a vortex
	crystal to give rise to a minimal surface, we need the nondegeneracy
	condition that allows the use of Implicit Function Theorem in our
	construction.  The connection is not clear for degenerate vortex crystals,
	which may or may not give rise to minimal surfaces.
\end{remark}

Unfortunately, we are aware of very few finite, nondegenerate, translating
vortex crystals; see Example~\ref{ex:translating}.  So the theorem above does
not bring us many new minimal surfaces.

Nevertheless, in Example~\ref{ex:translating}, we will obtain new minimal
surfaces from translating vortex crystals with imposed symmetry.  For a vortex
crystal, a \emph{circulation-preserving} (resp.\ \emph{-reversing}) symmetry is
a Euclidean isometry that maps vortices to identical (resp.\ opposite)
vortices.  When a symmetry group is imposed, the vortex crystal is said to be
\emph{nondegenerate} if the only perturbations that preserve the balance as
well as the symmetry are the trivial ones, namely Euclidean translations for
translating vortex crystals.  We have the following version of
Theorem~\ref{thm:translate} with imposed symmetry.  

\begin{theorem}[Translating vortex crystals with imposed symmetry]\label{thm:symmetric}
	Let $(p^\circ_k, \sigma_k)_{1 \le k \le n}$ be a normalized translating
	vortex crystal with $n$ vortices at $p^\circ_1$, \ldots, $p^\circ_n \in \C$.
	Let $G$ be a symmetry group of $(p^\circ_k, \sigma_k)$.  If the vortex
	crystal is nondegenerate with the symmetry group $G$ imposed, then the
	conclusion of Theorem~\ref{thm:translate} holds.  Moreover, the symmetry
	group $G$ induces a symmetry group of the resulting minimal surfaces.
\end{theorem}

In particular, a circulation-reversing reflection in the vortex crystal induces
an order-2 rotational symmetry around a straight line in the minimal surface,
and a circulation-preserving reflection induces a reflection symmetry for the
minimal surface.

\medskip

We did not manage to establish a similar connection between singly periodic
minimal surfaces and finite stationary vortex crystals.  See
Remark~\ref{rmk:stationary2} for detailed explanation.

\medskip

Periodic vortex crystals will give rise to doubly or triply periodic minimal
surfaces, as stated in the following theorems.

\begin{theorem}[Singly periodic vortex crystals]\label{thm:singly}
	Let $(p^\circ_k, \sigma_k)_{1 \le k \le n}$ be a nondegenerate vortex
	crystal with $n$ vortices at $p^\circ_1$, \ldots, $p^\circ_n \in \C / \langle
	1 \rangle$.  Then there exists a one-parameter family $(M_\varepsilon)_{0 <
	\varepsilon < \delta}$ of embedded minimal surfaces in $\R^3$ such that:
	\begin{enumerate}
		\item $M_\varepsilon$ admits translational symmetries along the vectors
			\[
				T_{0,\varepsilon} = 2\pi(\varepsilon\re\nu, \varepsilon\im\nu, 1)
				\quad\text{and}\quad
				T_{1,\varepsilon} = (\varepsilon^{-1}, 0, m\pi),
			\]
			where $m = \sum\sigma_k$, $\nu$ is related to the velocity $v$ of the
			translating vortex crystal by $\nu = -2\pi v$, and $\nu = 0$ if the vortex
			crystal is stationary.  So $M_\varepsilon$ is a doubly periodic minimal
			surface.

		\item The quotient $M_\varepsilon / \langle T_{0,\varepsilon},
			T_{1,\varepsilon}\rangle$ is of genus $n-1$ and has four Scherk ends
			(asymptotic to half-planes).

		\item The flux vector along any closed curve in $M_\varepsilon / \langle
			T_{0,\varepsilon}, T_{1, \varepsilon} \rangle$ has no vertical component.

		\item As $\varepsilon \to 0$, up to a translation, $M_\varepsilon/\langle
			T_{1,\varepsilon}\rangle$ converges to a helicoid of period $(0, 0,
			2\pi)$ in the neighborhood of $(p^\circ_i/\varepsilon, 0)$ for each $1
			\le i \le n$.  The helicoid is right-handed (resp.\ left-handed) if
			$\sigma_i = 1$ (resp.\ $-1$).
	\end{enumerate}
\end{theorem}

Recall that the flux vector along a closed curve $\Gamma$ is defined as the
integral of the conormal vector along $\Gamma$~\cite{perez1993}.  It can be
physically interpreted as the surface tension force along the
curve~\cite{perez2002}.

In Theorem~\ref{thm:singly}, we actually constructed a 3-parameter family of
minimal surfaces, depending on a real parameter $\varepsilon$ and a complex
parameter $\nu$.

\begin{theorem}[Doubly periodic vortex crystals]\label{thm:doubly}
	Let $(p^\circ_k, \sigma_k)_{1 \le k \le n}$ be a nondegenerate vortex
	crystal with $n$ vortices at $p^\circ_1$, \ldots, $p^\circ_n \in \T = \C /
	\langle 1, \tau \rangle$.  Assume that $m = \sum \sigma_k = 0$.  Then there
	exists a one-parameter family $(M_\varepsilon)_{0 < \varepsilon < \delta}$ of
	embedded minimal surfaces in $\R^3$ such that:
	\begin{enumerate}
		\item $M_\varepsilon$ admits a translational symmetry along the vectors
			\begin{align*}
				T_{0,\varepsilon} &= 2\pi(\varepsilon\re\nu, \varepsilon\im\nu, 1),\\
				T_{1,\varepsilon} &= (\varepsilon^{-1}, 0, \Psi_1(\varepsilon)),\\
				T_{2,\varepsilon} &= (\varepsilon^{-1}\re\tau, \varepsilon^{-1}\im\tau,
				\Psi_2(\varepsilon)),
			\end{align*}
			where $\nu$ is related to the velocity $v$ of the translating vortex
			crystal by $\nu = -2\pi v$, and $\nu = 0$ if the vortex crystal is
			stationary.  So $M_\varepsilon$ is a triply periodic minimal surface.

		\item The quotient $M_\varepsilon / \langle T_{0,\varepsilon},
			T_{1,\varepsilon}, T_{2,\varepsilon}\rangle$ is of genus $n+1$.

		\item The flux vector along any closed curve in $M_\varepsilon / \langle
			T_{0,\varepsilon}, T_{1, \varepsilon}, T_{2,\varepsilon} \rangle$ has no
			vertical component.

		\item As $\varepsilon \to 0$, up to a translation, $M_\varepsilon / \langle
			T_{1,\varepsilon}, T_{2,\varepsilon}\rangle$ converges to a helicoid of
			period $(0, 0, 2\pi)$ in the neighborhood of $(p^\circ_i/\varepsilon, 0)$
			for each $1 \le i \le n$.  The helicoid is right-handed (resp.\
			left-handed) if $\sigma_i = 1$ (resp.\ $-1$).  Moreover, we have
			\[
				\Psi_1(\varepsilon) \to -2\pi y \;\text{and}\;
				\Psi_2(\varepsilon) \to 2\pi x \;\text{as $\varepsilon \to 0$},
			\]
			where $(x, y) \in \R^2$ are defined by $\sum\sigma_kp_k = x + y \tau$.
	\end{enumerate}
	Conversely, for any such family of triply periodic minimal surfaces, in the
	limit $\varepsilon\to 0$, the positions and handednesses of the helicoids
	give a doubly periodic vortex crystal (possibly degenerate).
\end{theorem}

In Theorem~\ref{thm:doubly}, we actually constructed a 5-parameter family of
minimal surfaces, depending on a real parameter $\varepsilon$ and two complex
parameter $\nu$ and $\tau$.

\begin{remark}
	In the Theorems, the surface is rotated into a position so that the flux
	vectors are horizontal.  We find this choice best to reveal the connection to
	vortex crystals.  The price is that $\Psi_1$ and $\Psi_2$ are left to vary
	with $\varepsilon$.  One could also rotate the surface to fix $\Psi_1 =
	\Psi_2 \equiv 0$.  Then the flux vectors are not horizontal, and the
	connection to vortex crystals is less direct.
\end{remark}

\section{Examples}\label{sec:example}

The study on vortex crystals traces back to about 150 years ago, and has
accumulated plenty of examples, many of which are binary hence would imply
minimal surfaces.  A nice survey of these examples is provided by Aref et
al.~\cite{aref2003}.  Here, we examine the known vortex crystals and their
corresponding minimal surfaces.  Occasionally, we also have minimal surfaces
that lead to new vortex crystals.

\begin{example}[Linear configuration]
	Traizet and Weber~\cite{traizet2005} considered configurations of helicoids
	along a straight line.  In particular, when $p^\circ_i$, $1 \le i \le n$, are
	the roots of $H_n$, the Hermite polynomial of degree $n$, and $\sigma_i =
	-1$ for all $i$, then the configuration $(p^\circ_i, \sigma_i)_{1 \le i \le
	n}$ is balanced and nondegenerate.  In fluid dynamics, the corresponding
	rotating vortex crystal was first found by Stieltjes in
	1885~\cite{stieltjes1885}, and has been rediscovered and revisited many
	times~\cites{szego1959, marden1949, eshelby1951}.

	Traizet and Weber~\cite{traizet2005} also considered a configuration with $n
	= 2m+1$ helicoids, $m+1$ of which lie at the roots of $H_{m+1}$ and have
	positive handedness, and the remaining $m$ lie at the roots $H_m$ and have
	negative handedness.  This configuration is proved to be balanced and
	nondegenerate.  We are not aware of any discussion of corresponding vortex
	crystals in the fluid dynamics community.  So this is a new example of vortex
	crystal inspired by the minimal surfaces. \qed
\end{example}

\begin{figure}
	\includegraphics[height=0.3\textwidth]{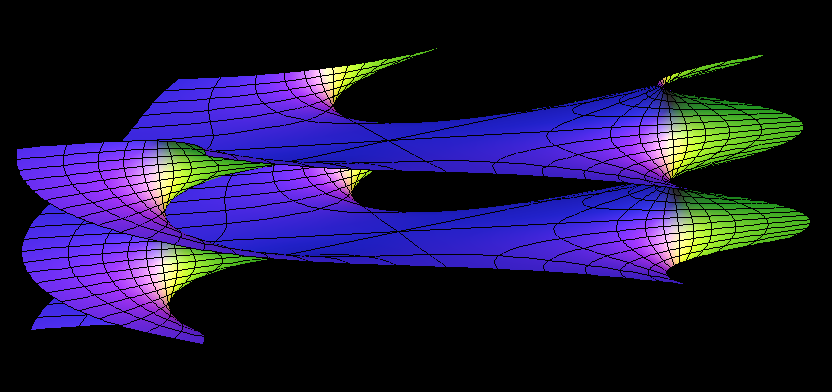}
	\includegraphics[height=0.3\textwidth]{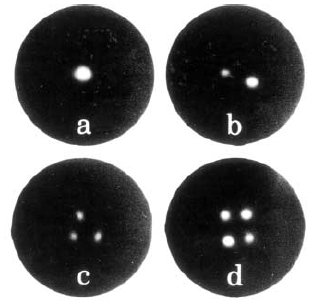}
	\caption{
		Left: A twisted Karcher--Scherk tower of six wings near the limit
		consisting of three helicoids at the vertices of an equiangular triangle
		(source: 3D-XplorMath Minimal Surface Gallery).  Right: Thomson's vortex
		polygons seen in superfluid Helium, reprinted with permission from
		\cite{yarmchuk1979}.  \label{fig:Scherk}
	}
\end{figure}

\begin{example}[Polygonal configuration]
	A Karcher--Scherk tower with $2n$ wings can be twisted into a configuration
	of $n$ negatively handed helicoids lying at the vertices of a regular
	polygon; see Figure~\ref{fig:Scherk} and~\cite{freese2021}*{Prop.~8.11}.  The
	Fischer--Koch surfaces can be twisted into a similar configuration, only with
	an extra helicoid, positively or negatively handed, at the center of the
	polygon; see~\cite{freese2021}*{Prop.~8.13}.

	The corresponding rotating vortex crystals were first investigated by Thomson
	in 1883~\cite{thomson1883}.  In particular, he famously proved that identical
	vortices at the vertices of a regular $n$-gon is linearly stable if $n \le
	6$, and linearly unstable if $n \ge 8$, and ``neutrally stable'' if $n = 7$.
	The proof was improved and modified many times~\cites{havelock1931,
	dritschel1985, aref1995}.

	For our minimal surfaces, this stability analysis means that the Helicoid
	limit of twisted Scherk surface is nondegenerate for $n \ne 7$, for which
	Theorem~\ref{thm:translate} applies.  The $n = 7$ case is degenerate, but
	becomes nondegenerate if we impose the dihedral symmetry of the
	heptagon~\cite{freese2021}, so Theorem~\ref{thm:symmetric} applies.  Even
	without imposed symmetry, it was proved that the nonlinear stability still
	holds for $n=7$~\cite{kurakin2002}, so a more elaborated version of the
	implicit function theorem might apply.  \qed
\end{example}

\begin{example}[Nested polygonal configurations]
	Fluid dynamists also investigated vortex crystals where vortices lie on the
	vertices of several concentric regular polygons~\cites{havelock1931,
	campbell1978, aref2005}.  In some cases, e.g.\ when the vortices lie on two
	concentric polygons, with or without an extra vortex at the center, an
	algebraic approach is possible.  However, we are not aware of any systematic
	investigation on such configurations.

	This line of research overlaps with~\cite{freese2021}, where the second named
	author considered configurations with dihedral symmetry, and helicoids all
	lie on the symmetry lines (including the center).  Each dihedral
	configuration of helicoids can be seen as corresponding to a nested polygonal
	vortex crystal.

	For instance, the Callahan-Hoffman-Meeks (CHM) surface~\cite{callahan1989} of
	genus $2k+1$, $k \ge 1$, can be twisted into a family of minimal surfaces
	invariant under screw-motions~\cite{callahan1990}.  Numerical computation
	suggests a smooth deformation with the angle $\theta$ of the screw-symmetry
	varying within the open interval $(-\frac{2\pi}{k+1},
	\frac{2\pi}{k+1})$~\cite{callahan1993}.  As $\theta$ approaches $\pm
	2\pi/(k+1)$, the surface seems to diverge to a degenerate
	limit~\cite{callahan1993}.  Computer images suggest that this is a helicoid
	limit corresponding to a vortex crystal with positive vortices forming a
	regular $(k+1)$-gon, and negative vortices forming a concentric regular
	$(k+1)$-gon.  One may determine from the balance
	equations~\eqref{eq:balancevortex} that the ratio $r$ of the circumradii of
	the polygons is one of the two solutions (inverse to each other) to
	$k\frac{1-r^k}{1+r^k} = \frac{1+r^2}{1-r^2}$.  See~\cite{freese2021}*{Figure
	1} for an image with $k=1$.  We are not aware of any discussion of these
	vortex crystals in the fluid dynamics community, so they might be new. \qed
\end{example}

\begin{example}[Numerical examples]
	Campbell and Ziff~\cites{campbell1978, campbell1979} have obtained numerical
	examples of vortex crystals, and claimed to have found linearly stable
	configurations with up to $30$ identical vortices.  Their 1978 report is
	often referred to as the \emph{Los Alamos Catalog}.  Many of their examples
	were experimentally observed in superfluid Helium~\cite{yarmchuk1979}; see
	Figure~\ref{fig:Scherk}, right.  By our connection, they all correspond to
	screw-motion invariant minimal surfaces.

	In many of their examples, the vortices seem to lie on concentric rings, but
	this impression is not precise.  Rather, all their examples admit an axis of
	symmetry.  Asymmetric examples were not found until~\cite{aref1998}. \qed
\end{example}

\begin{figure}
	\includegraphics[width=0.6\textwidth]{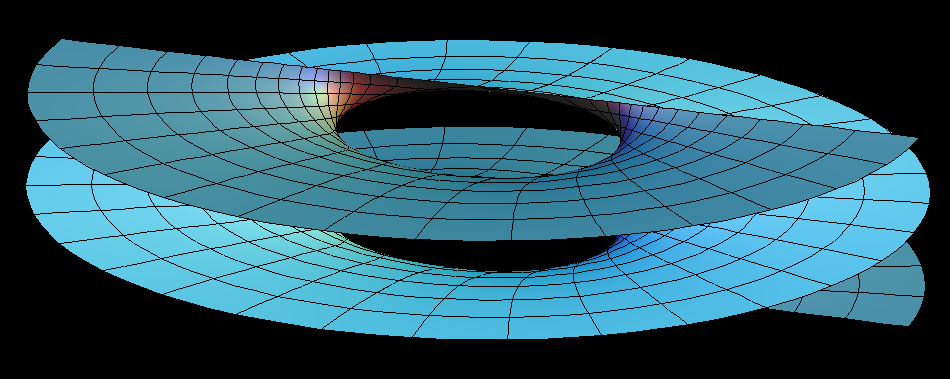}
	\includegraphics[width=0.6\textwidth]{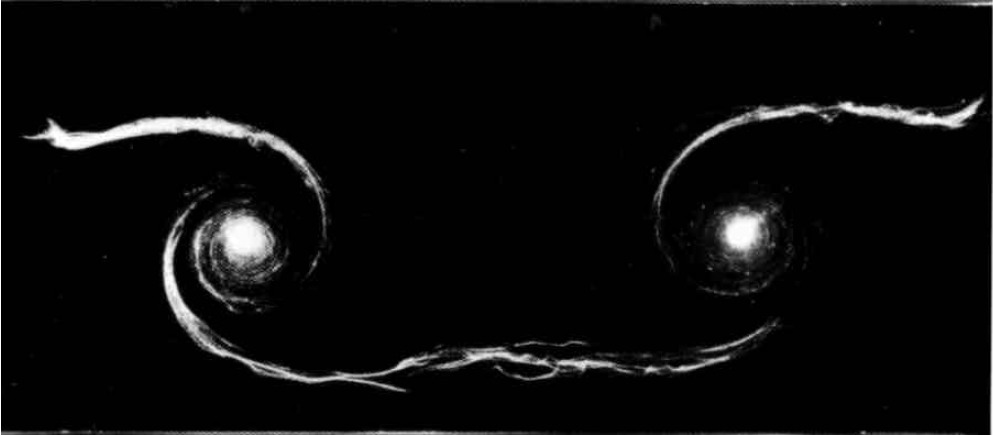}
	\caption{
		Top: A Riemann minimal example near the limit consisting of a pair of
		opposite helicoids (source: 3D-XplorMath Minimal Surface Gallery).  Bottom:
		A pair of opposite vortices seen in a cross-section of the vortex sheet
		behind a wing. (source: H.\ Bippes~\cite{vandyke1982}*{p.\ 50}).
		\label{fig:Riemann}
	}
\end{figure}

\begin{figure}
	\includegraphics[width=0.95\textwidth]{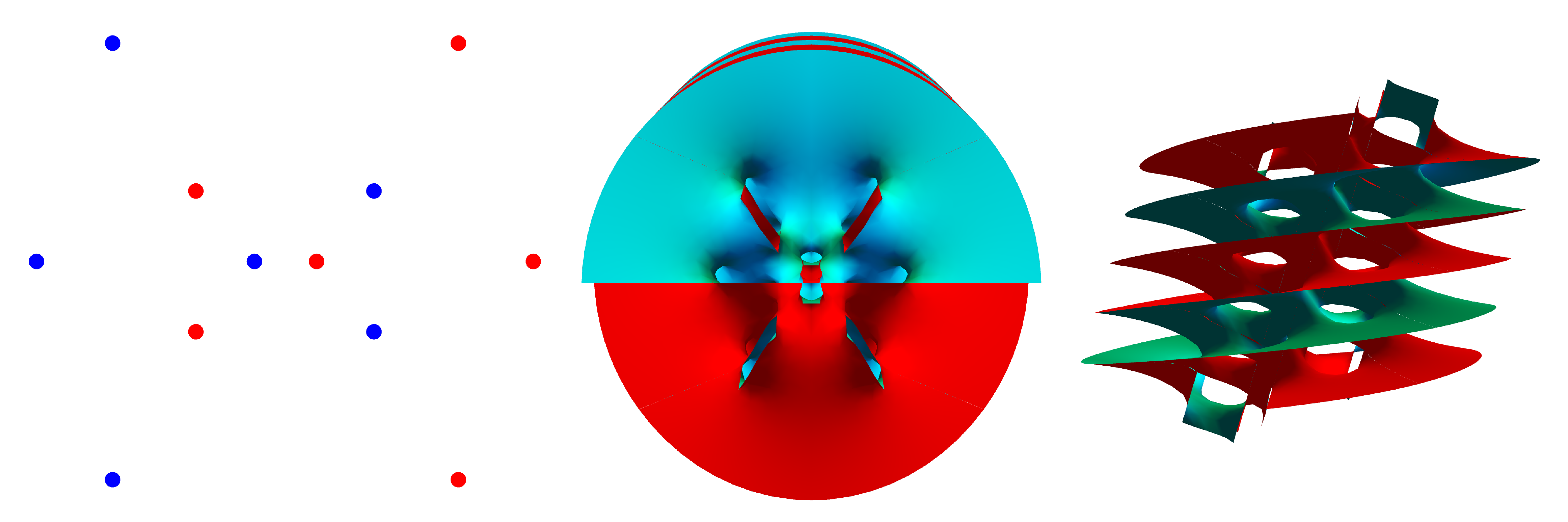}
	\caption{
		Left: An Adler--Moser translating vortex crystal with $n/2=6$.  Middle: Top
		view of a singly periodic minimal surface arising from the corresponding
		helicoid configuration.  Right: Side view of the same minimal surface.
		\label{fig:translate}
	}
\end{figure}

\begin{example}[Translating vortex crystal]\label{ex:translating}
	A translating vortex crystal must consist of an even number $n$ of vortices,
	with $n/2$ positive vortices and $n/2$ negative ones.  The simplest case is a
	pair of opposite vortices, corresponding to Riemann minimal examples; see
	Figure~\ref{fig:Riemann}.  Surprisingly, there is no solution for $n/2 = 2$.
	Some examples for $n/2 = 3$ and $n/2 = 6$ can be found in~\cite{kadtke1987}.  

	More generally, when $n/2=j(j+1)/2$ is a triangular number, binary
	translating vortex crystals have been found~\cites{bartman1984, campbell1987,
	clarkson2009} with positive vortices at the roots of a $j$-th Adler--Moser
	polynomial $\Theta_j$, and negative vortices at the roots of the
	corresponding modified Adler--Moser polynomial $\tilde \Theta_j$.

	The Adler--Moser polynomial $\Theta_j(z)$ actually depends on $m$ complex
	parameters $\kappa_1 = z$, $\kappa_2$, \ldots, $\kappa_j$.  Changing these
	parameters preserves the balance.  In fact, these perturbations are linearly
	independent and span the kernel of the Jacobian~\cite{liu2020}.  As a
	consequence, the Adler--Moser translating configurations are degenerate
	except for the trivial case $n=2$.

	However, under the assumption that $\Theta_j$ has only simple roots, these
	configurations are proved to be nondegenerate~\cite{liu2020} if we impose (up
	to Euclidean motions) a reflection symmetry in the real axis that preserves
	circulations and a reflection symmetry in the imaginary axis that reverses
	circulations.  In fact, this symmetry is realized by a unique choice of
	parameters $\kappa_2$, \ldots, $\kappa_j$.  With this choice, the assumption
	that $\Theta_j$ has only simple roots is verified for $j \le
	34$~\cite{liu2020}.

	By Theorem~\ref{thm:symmetric}, these symmetric translating Adler--Moser
	configurations give rise to singly periodic minimal surfaces.  The reflection
	in the real axis becomes a rotational symmetry in the $x$-axis, and the
	reflection in the imaginary axis becomes a reflection in the $yz$-plane.  See
	Figure~\ref{fig:translate} for an example with $n/2=6$.  \qed
\end{example}

\begin{figure}
	\includegraphics[width=0.6\textwidth]{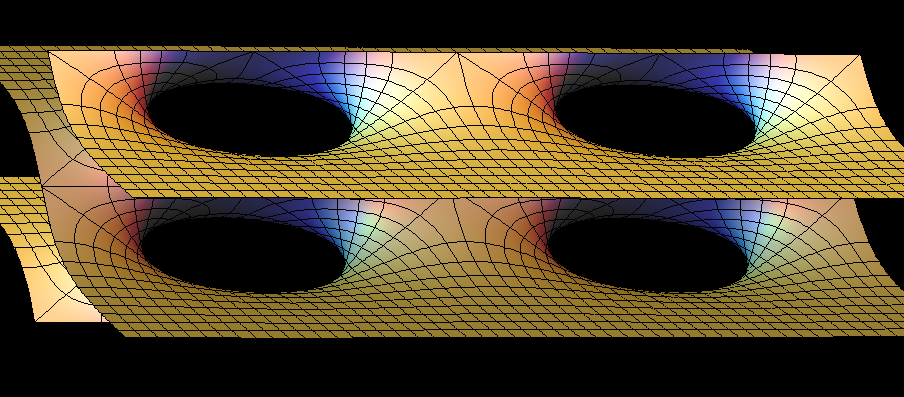}
	\includegraphics[width=0.6\textwidth]{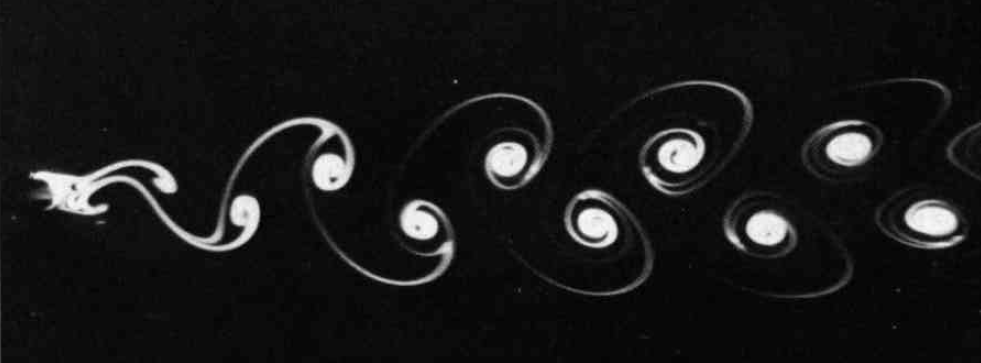}
	\caption{
		Top: A doubly periodic Karcher--Meeks--Rosenberg surface near the limit
		consisting of a row of alternating helicoids (source: 3D-XplorMath Minimal
		Surface Gallery).  Bottom: A von-K\'arm\'an vortex street.  (source:
		Sadatoshi Taneda~\cite{vandyke1982}*{p.\ 57}).  \label{fig:KMR}
	}
\end{figure}

\begin{example}[Singly periodic vortex crystals]
	The doubly periodic Scherk surface can be deformed to a periodic helicoid
	limit.  It corresponds to a singly periodic vortex crystal with a single
	vortex in the period.

	The famous vortex street of von K\'arm\'an~\cite{vonkarman1912}, and more
	general cases of Dolaptschiew and Maue~\cite{maue1940}, are the only singly
	periodic vortex crystals with two (opposite) vortices in the period.  They
	correspond to the helicoid limits of Karcher--Meeks--Rosenberg
	surfaces~\cites{karcher1988, meeks1988}.  See Figure~\ref{fig:KMR}. \qed
\end{example}

\begin{figure}
	\includegraphics[height=0.4\textwidth]{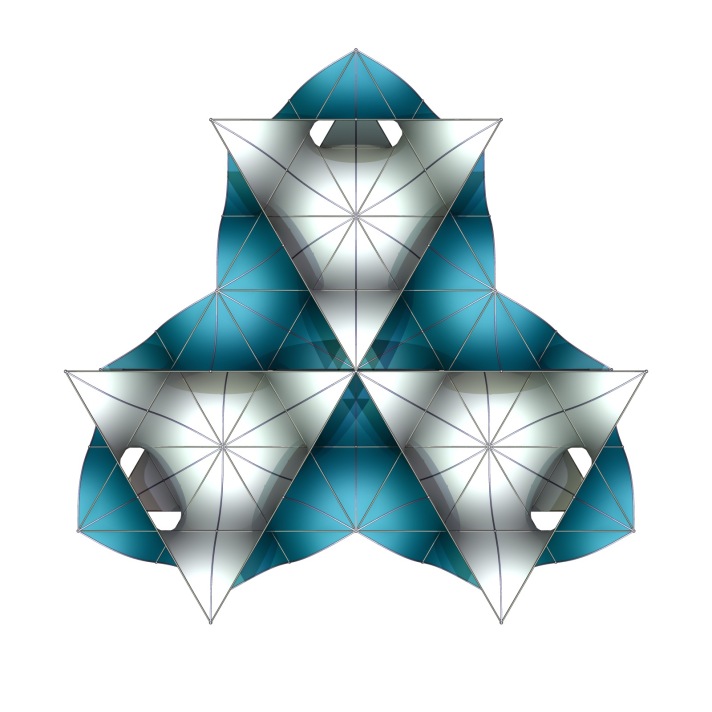}
	\includegraphics[height=0.4\textwidth]{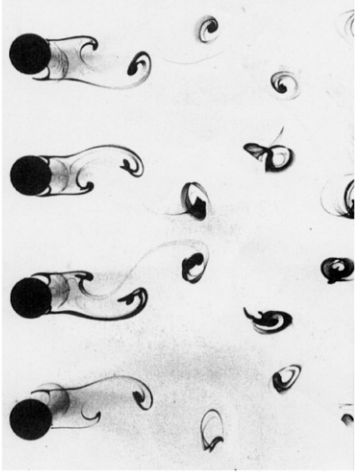}
	\caption{
		Left: A triply periodic rPD surface near the limit consisting of helicoids
		arranged in a hexagonal lattice (source: Matthias Weber).  Right: A doubly
		periodic vortex crystal seen in the wake behind a row of cylinders.
		(source: Toshio Kobayashi~\cite{kobayashi1984}*{p.\ 43}).  \label{fig:rPD}
	}
\end{figure}

\begin{example}[Doubly periodic vortex crystals]
	If a doubly periodic vortex crystal has two vortices in the period, they must
	be of opposite circulations.  Such a configuration is generically
	nondegenerate.  They give rise to triply periodic minimal surfaces of genus
	3 (TPMSg3); see Figure~\ref{fig:rPD} for an example.  In fact, our
	construction will assume an orientation-reversing translation, hence the
	produced examples must all belong to the 5-parameter family of TPMSg3s
	constructed by Meeks~\cite{meeks1990}. \qed
\end{example}

\section{Sketched construction}\label{sec:construct}

In~\cites{traizet2005, freese2021}, the main technical issue was the
multivaluedness of the Weierstrass data.  For translating invariant minimal
surfaces, the Weierstrass data is single-valued, so the construction is much
easier.  The construction below is adapted from~\cite{traizet2008}.  This
approach has been repeated many times in the literature~\cites{traizet2008,
traizet2008b, chen2021} and is only simpler in our context; hence we will only
present a sketch.

\subsection{Weierstrass data}

Let $\Sigma_+$ be
\begin{itemize}
	\item the Riemann sphere $\hat\C = \C \cup \{\infty\}$ if the vortex crystal is finite;

	\item the annuli $\C / \langle 1 \rangle$ if the vortex crystal is singly periodic;

	\item the torus $ \T_+ = \C / \langle 1, \tau_+ \rangle$ if the vortex crystal is
		doubly periodic.
\end{itemize}
Moreover, let $\Sigma_-$ be
\begin{itemize}
	\item Another copy of $\Sigma_+$ if the vortex crystal is finite or singly periodic;

	\item the torus $ \T_- = \C / \langle 1, \tau_- \rangle$ if the vortex crystal is
		doubly periodic.
\end{itemize}
Consider $n$ points $p=(p_k)_{1 \le k \le n}$ in $\Sigma_+$ and $n$ points
$q = (q_k)_{1 \le k \le n}$ in $\Sigma_-$. 

The node-opening is parameterized by $n$ complex numbers $t = (t_k)_{1 \le k
\le n}$.  If $t=0$, we identify $p_k \in \C_+$ and $q_k \in \C_-$ to form a
node.  The resulting singular Riemann surface with nodes is denoted $\Sigma_0$.
If $t \ne 0$, we open the nodes as follows.  Let $z_\pm$ be the standard
coordinates of $\C_\pm$.  Consider local coordinates $w_k^+ = z_+ - p_k$ in the
neighborhood of $p_k$ and $w_k^- = z_- - q_k$ in the neighborhood of $q_k$.
Fix a small $\delta>0$ such that the disks $|w_k^\pm| < \delta$ are all
disjoint.  Then for every $1 \le k \le n$, we remove the disk $|w_{k,\pm}| <
|t_k|/\delta$, and identify the annuli
\[
	|t_k|/\delta \le |w_k^+| \le \delta \qquad \text{and} \qquad
	|t_k|/\delta \le |w_k^-| \le \delta
\]
by
\[
	w_k^+ w_k^- = t_k.
\]
This produces a Riemann surface that we denote $\Sigma_t$.

We construct the minimal surface using the Weierstrass parameterization
\[
	\Sigma_t \ni z \mapsto \re\int_{z_0}^z \big(\phi_1, \phi_2, \phi_3)
\]
where $\phi_1$, $\phi_2$, and $\phi_3$ are meromorphic forms on $\Sigma_t$
satisfying the conformality condition
\begin{equation}\label{eq:conformal}
	Q = \phi_1^2 + \phi_2^2 + \phi_3^2 = 0.
\end{equation}

Then the flux vector along a closed curve $\Gamma$ is given by~\cite{perez1993}

\begin{align*}
  \im \int_\Gamma \left(\phi_1, \phi_2,\phi_3\right).
\end{align*}

\subsection{Equations}

We define
\[
	\Omega_\pm = \{ z \in \Sigma_\pm \mid |w_k^\pm| > \delta\, \text{for all}\, 1 \le k \le n \}.
\]
Let $A_k$ be an anticlockwise circle in $\Omega_+$ around $p_k$ and $A'_k$ be
an anticlockwise circle in $\Omega_-$ around $q_k$.  Note that $A_k$ is
homologous in $\Sigma_t$ to $-A'_k$.  Let $B_k$ be a cycle in
$\Sigma_\varepsilon$ which goes ``half-way up'' from $\Sigma_+$ to $\Sigma_-$
through the helicoind near $p_1$ then ``half-way down'' through the helicoid
near $p_k$, as in \cite{freese2021}.  We need to solve the period problems
\begin{align}
	\re \int_{A_k} (\phi_1, \phi_2, \phi_3) &= 2\pi\sigma_k (\varepsilon \re \nu, \varepsilon \im \nu, 1),\nonumber\\
	\re \int_{B_k} (\phi_1, \phi_2, \phi_3) &= (0, 0, 0)\label{eq:periodBk}.
\end{align}
We close the $A$-period by defining $\phi$'s as the unique meromorphic forms
satisfying
\[
	\int_{A_k} (\tphi_1, \tphi_2, \tphi_3) =
	2\pi\ii (\alpha_k - \ii \sigma_k \varepsilon^2 \re\nu, \beta_k - \ii \sigma_k \varepsilon^2 \im\nu, \gamma_k - \ii\sigma_k\varepsilon),
\]
where $\tphi_i := \varepsilon\phi_i$ are the rescaled Weierstrass data.
Depending on the type of the vortex crystal, we also require the following
\begin{itemize}[leftmargin=*]
	\item if the vortex crystal is finite and translating, we want $\tphi_1$ and
		$\tphi_2$ to have double poles at $\infty_\pm$.  Up to rotations and
		scalings, we assume that
		\begin{equation}\label{eq:assumefinite}
			\tphi_1 \sim dz_\pm + \mathcal{O}(z^{-2}_\pm)dz_\pm
			\quad \text{and} \quad
			\tphi_2 \sim \mp \ii dz_\pm + \mathcal{O}(z^{-2}_\pm)dz_\pm
			\quad \text{at $\infty_\pm$}.
		\end{equation}
		On the other hand, since the minimal surfaces have planar ends, $\tphi_3$
		must be holomorphic at $\infty_\pm$.  Consequently, we must have $m = \sum
		\sigma_k = 0$.


	\item if the vortex crystal is singly periodic, let $A_\pm \subset
		\Sigma_\pm$ be the segments $\{t \mp K \ii \mid 0 \le t \le 1\}$, where $K
		> |\im p_k|$ and $K > |\im q_k|$ for all $k$.  We want that
		\[
			\re \int_{A_\pm} (\tphi_1, \tphi_2, \tphi_3) =
			(1, 0, \varepsilon m \pi).
		\]
		So we require that
		\[
			\int_{A_\pm} (\tphi_1, \tphi_2, \tphi_3) =
			(1 + \ii \alpha_\pm, \ii \beta_\pm, \varepsilon m \pi + \ii \gamma_\pm).
		\]
		The flux vectors $(\alpha_\pm, \beta_\pm, \gamma_\pm)$ (or rather
		their inverse) can be physically interpreted as the surface tension forces
		along the Scherk ends.  Up to a rotation around the x-axis (corresponding
		to the real axis), we may assume that $\gamma_+ \equiv 0$.

	\item if the vortex crystal is doubly periodic, let $A_\pm \subset
		\Omega_\pm$ be curves homologous in $\Sigma_\pm$ to the segment from $0$ to
		$1$, and $B_\pm \subset \Omega_\pm$ be curves homologous to the segment
		from $0$ to $\tau_\pm$.  We want that
		\begin{align}
			\re \int_{A_\pm} (\tphi_1, \tphi_2, \tphi_3) &=
			(1, 0, \varepsilon\Psi_1),\nonumber\\
			\re \int_{B_\pm} (\tphi_1, \tphi_2, \tphi_3) &=
			(\re \tau, \im \tau, \varepsilon\Psi_2).\label{eq:periodB}
		\end{align}
		So we require that
		\[
			\int_{A_\pm} (\tphi_1, \tphi_2, \tphi_3) =
			(1 + \ii \alpha_\pm, \ii \beta_\pm, \varepsilon\Psi_1 + \ii \gamma_\pm).
		\]
		Again, the flux vectors can be physically interpreted as the surface
		tension forces.  Up to a Euclidean rotation, we may assume that the periods
		of $\tphi_3$ over $A_+$ and $B_+$ are real.  That is
		\[
			\gamma_+ = \im \int_{A_+} \tphi_3 = 0
			\quad\text{and}\quad 
 			\im \int_{B_+} \tphi_3 = 0.
		\]
\end{itemize}
In any of these cases, the $\phi$'s are uniquely determined by the requirements
above.

\medskip

Write $\tQ:=\varepsilon^2 Q$.  The conformality condition~\eqref{eq:conformal}
is equivalent to
\begin{align}
	\mathcal{E}_k &:= \int_{A_k} \frac{\tQ w_k^+}{dz_+} = 0, & 1 \le k \le n, \label{eq:conformalAk}\\
	\mathcal{F}_k &:= \int_{A_k} \frac{\tQ}{dz_+} = 0, & 1 \le k \le n, \label{eq:conformalAkp}\\
	\mathcal{F}'_k &:= \int_{A'_k} \frac{\tQ}{dz_-} = 0, & 1 \le k \le n, \label{eq:conformalAkq}
\end{align}
and, if the vortex crystal is periodic,
\begin{equation}\label{eq:conformalA}
	\int_{A_\pm} \frac{\tQ}{dz_\pm} = 0.
\end{equation}
Note that \eqref{eq:conformalAkp} are not independent.  One dependence comes
from the residue theorem, namely that
\begin{equation}\label{eq:depend1}
	\sum_{k=1}^n \mathcal{F}_k = 0.
\end{equation}

\subsection{Solutions}

To facilitate the solution, we will construct minimal surfaces with an
orientation-reversing translational symmetry $R_\varepsilon$ such that
$R_\varepsilon^2 = T_{0,\varepsilon}$.  We want $R_\varepsilon$ to correspond
to the symmetry
\[
	\iota: \Sigma_+ \ni z \mapsto \overline{z} \in \Sigma_-.
\]
More specifically, we want that
\[
	\iota^*(\phi_1, \phi_2, \phi_3) = (\overline{\phi_1}, \overline{\phi_2}, \overline{\phi_3}).
\]
This can be achieved by assuming that
\[
	q_k \equiv \overline{p}_k, \quad
	t_k \in \R \quad\text{and}\quad
	(\alpha_k, \beta_k, \gamma_k) \equiv 0.
\]
In the doubly periodic case, we also assume that $\tau_- \equiv
\overline{\tau_+}$.  As a consequence, the $B$-period
problem~\eqref{eq:periodBk} are automatically solved.  Indeed, we can choose
$B$-curves such that $B_k + \iota(B_k)$ is homologous to $\sigma_1 A_1 -
\sigma_k A_k$, for which all $\phi$ periods are pure imaginary;
see~\cite{freese2021}.  In the periodic cases, since $\iota(A_+) = A_-$, we
also have $\gamma_- = \gamma_+ \equiv 0$.  In the doubly periodic case, since
$\iota(B_+) = B_-$, the period of $\tphi_3$ over $B_-$ must also be real.  Note
that the $A$-curves and $B$-curves form a homology basis.  Now that the
vertical components of the flux vectors vanish along all these curves, they
must vanish along any closed curve.

Moreover, since $\iota^*(Q) = \overline{Q}$ and $\iota(A_k) = A_k = -A'_k$, we
have
\[
	\mathcal{E}_k = \int_{A_k} \frac{\tQ w_k^+}{dz_+}
	= \int_{\iota(A_k)} \iota^* \Big( \frac{\tQ w_k^+}{dz_+}\Big)
	= -\int_{A'_k} \frac{\overline\tQ w_k^-}{dz_-}
	= -\overline{\mathcal{E}_k},
\]
and
\[
	\mathcal{F}_k = \int_{A_k} \frac{\tQ}{dz_+}
	= \int_{\iota(A_k)} \iota^*\Big(\frac{\tQ}{dz_+}\Big)
	= -\int_{A'_k} \frac{\overline{\tQ}}{dz_-}
	= -\overline{\mathcal{F}'_k}.
\]
This means that $\re \mathcal{E}_k = 0$ and~\eqref{eq:conformalAkq} is
automatically solved if~\eqref{eq:conformalAkp} is solved.

\begin{remark}\label{rmk:iota1}
	Without the symmetry $\iota$, we can use the Implicit Function Theorem to
	prove the existence of $q_k$, $\im t_k$, $\alpha_k$, $\beta_k$, $\gamma_k$,
	and $\tau_-$ that solve the $B$-period problem and the conformality
	equation~\eqref{eq:conformalAkq}.
\end{remark}

Similar arguments as in~\cite{traizet2008}, using the Implicit Function
Theorem, show that

\begin{proposition}\label{prop:alphabetat}
	For $\varepsilon$ in a neighborhood of $0$ and $p_k$ in a neighborhood of its
	central values $p^\circ_k$, there exist unique values for parameters $t_k$,
	$\alpha_\pm$, $\beta_\pm$, and $\tau_+$, depending smoothly on $\varepsilon$
	and $p_k$, such that the imaginary part of~\eqref{eq:conformalAk}, as
	well as~\eqref{eq:conformalA} in the periodic cases,
	and~\eqref{eq:periodB} in the doubly periodic case, are solved.  At
	$\varepsilon=0$, we have $t_k = 0$ in all cases, $\alpha_\pm = 0$, $\beta_\pm
	= \mp 1$ in the periodic cases, and $\tau_+ = \tau$ in the doubly periodic
	case, where $\tau$ is the torus parameter for the given doubly periodic
	vortex crystal.
\end{proposition}

A sketched proof for the proposition is delayed to Appendix~\ref{app:alphabetat}.

At $\varepsilon = 0$, we have
\[
	\tphi_1 = dz_\pm \quad\text{and}\quad \tphi_2 = \mp \ii dz_\pm \quad\text{on}\quad\Sigma_\pm.
\]
On $\Sigma_+$, $\phi_3 = \tphi_3 / \varepsilon$ extends smoothly to
$\varepsilon = 0$ with the explicit form
\begin{equation}\label{eq:phi3}
	\sum -\ii\sigma_k\Upsilon(z-p_k)dz
\end{equation}
where \footnote{Note that, in the doubly periodic case, $\Upsilon$ is not
meromorphic, but $\phi_3$ is.}
\[
	\Upsilon(z) :=
	\begin{dcases}
		1/z &
		\text{if $(p_k,\sigma_k)$ is finite};\\
		\pi \cot(\pi z) &
		\text{if $(p_k,\sigma_k)$ is singly periodic};\\
		\big(\zeta(z; \tau) - \xi(z; \tau)\big) &
		\text{if $(p_k,\sigma_k)$ is doubly periodic.}
	\end{dcases}
\]
One then verifies that, at $\varepsilon = 0$, we have indeed
\[
	\int_{A_\pm} \phi_3 = m\pi \in \R
\]
if $(p_k, \sigma_k)$ is singly periodic and (cf.~\cite{traizet2008}*{\S~4.3.1}
and \cite{chen2021}*{\S~5})
\[
	\int_{A_\pm} \phi_3 = -2\pi y \in \R \quad\text{and}\quad
	\int_{B_\pm} \phi_3 = 2\pi x \in \R
\]
if $(p_k,\sigma_k)$ is doubly periodic, as we have assumed.  Here $(x, y) \in
\R^2$ are defined by $\sum\sigma_kp_k = x + y \tau$.

Then, as $\varepsilon \to 0$, $\mathcal{F}_{j,+} / \varepsilon^2$ converges
smoothly to the value
\begin{align*}
	\frac{\partial\mathcal F_{j,+}}{\partial \varepsilon^2}\Big|_{\varepsilon= 0}
	&= \int_{A_j} \Big(\frac{2\tphi_1}{dz_+} \frac{\partial \tphi_1}{\partial \varepsilon^2} + \frac{2\tphi_2}{dz_+} \frac{\partial \tphi_2}{\partial \varepsilon^2} \Big)+ 2\pi i \res_{p_j} \frac{\phi_3^2}{dz_+}\\
	&= \frac{\partial }{\partial \varepsilon^2}\int_{A_j} (2 \tphi_1 - 2i\tphi_2 )
	- 2\pi i \res_{p_j}\Big(\sum_k \sigma_k \Upsilon(z-p_k)\Big)^2\\
	&=4\pi\sigma_j\re \nu - 4\pi\ii\sigma_j \im \nu
	-4\pi\ii\sigma_j \sum_{k \ne j} \sigma_k \Upsilon(p_j-p_k)\\ &
	= 4\pi\sigma_j \Big[ \overline{\nu}
	-\ii \sum_{k \ne j} \sigma_k \Upsilon(p_j-p_k)\Big] = 8\pi^2 F_j,
\end{align*}
which vanishes if and only if $(p_k)$ and $(\sigma_k)$ are the positions and
circulations of a binary vortex crystal with velocity $v = -\nu / 2\pi$.  This
proves that the helicoid limits are necessarily vortex crystals.

Recall that $\mathcal{F}_k$ are related by~\eqref{eq:depend1}.  Therefore, if
the vortex crystal is not a finite stationary one, and is nondegenerate
(possibly with imposed symmetry), we may apply the Implicit Function Theorem to
prove the following proposition that concludes the construction of a family of
immersed surfaces.

\begin{proposition}
	If $(p^\circ_k)$ and $(\sigma_k)$ are the positions and circulations of a
	binary vortex crystal that is nondegenerate (possibly with imposed
	symmetry).  Then for $\varepsilon$ in a neighborhood of $0$, there exist
	unique values for $p_k$, depending smoothly on $\varepsilon$, such that
	$p_k(0) = p^\circ_k$ and the conformality condition~\eqref{eq:conformalAkp}
	is solved.
\end{proposition}

Finally, a similar argument as in~\cites{traizet2005,freese2021} proves that the
constructed surfaces are all embedded for $\varepsilon$ sufficiently small.

\begin{remark}
	The same construction can be carried out without the symmetry $\iota$ (see
	Remark~\ref{rmk:iota1}) and the conclusion is the same.  By uniqueness of
	the implicit functions, this implies that all minimal surfaces sufficiently
	close to a balanced and nondegenerate configuration of helicoids admit an
	orientation-reversing symmetry.
\end{remark}

\begin{remark}\label{rmk:stationary2}
	For finite stationary vortex crystals, we may fix two vortices to quotient
	out Euclidean similarities, leaving $n-2$ free complex parameters.  But in
	the last step of the construction, we have $n$ complex equations $F_k=0$ to
	solve.  The relation~\eqref{eq:depend1} eliminates one complex equation,
	Riemann Bilinear relation can eliminate a real equation
	(see~\cites{traizet2002, traizet2008b}), but we still have one real equation too
	many.  Hence the construction does not work.  This is compatible with
	Traizet's observation that catenoids cannot be glued into a single periodic
	minimal surface with vertical periods~\cite{traizet2002b}.
\end{remark}

\begin{remark}
	In~\cite{traizet2015}, the height differential of the minimal surface
	corresponds to the velocity field of the fluid.  A similar correspondence can
	be seen in our construction by noticing from~\eqref{eq:phi3} that, at
	$\varepsilon=0$, $\phi_3 = 2\pi \overline{u} dz$, where $u$ is the flow
	velocity field generated by the vortex crystal (away from the vortices).
\end{remark}

\appendix

\section{Proof of proposition~\ref{prop:alphabetat}}\label{app:alphabetat}

Similar arguments as in~\cite{traizet2008} apply here.

When $\varepsilon = 0$, all the $A_k$-periods of $\tphi$'s vanish, and thus
$\tphi_3$ converges to $0$.  If the vortex crystal if finite, $\tphi_1$ and
$\tphi_2$ converge to holomorphic forms in $\C$, and their limits are
determined by their behavior at $\infty$.  In view of the
assumptions~\eqref{eq:assumefinite}, we have
\[
	\tphi_1 \to dz_\pm \quad \text{and} \quad \tphi_2 \to
	\mp\ii dz_\pm 
\]
in $\Sigma_\pm$ as $\varepsilon \to 0$, and one easily verifies that
$\tQ$ converges to $0$.

If the vortex crystal is periodic, $\tphi_1$ and
$\tphi_2$ converge to holomorphic forms in $\Sigma_\pm$ and are
determined by their $A_\pm$ periods.  More specifically, at $\varepsilon =
0$, we have
\[
  \tphi_1 \to (1 + i\alpha_\pm)dz_\pm
  \quad \text{and} \quad
  \tphi_2 \to \ii\beta_\pm dz_\pm
\]
in $\Sigma_\pm$ as $\varepsilon \to 0$.  Whence we have,
\[
  \tQ = (\tphi_1)^2 + (\tphi_2)^2
  \to (1 + 2 i \alpha_\pm - \alpha_\pm^2 - \beta_\pm^2)dz_\pm^2.
\]
In order for this to vanish, we need $\alpha_\pm = 0$ and $\beta_\pm = \mp
1$.  The sign of $\beta_\pm$ is chosen so the surface has the desired
orientation.  Hence, again, we have the limit
\[
	\tphi_1 \to dz_\pm \quad \text{and} \quad \tphi_2 \to \mp\ii dz_\pm
\]
in $\Sigma_\pm$ as $\varepsilon \to 0$.

The remaining proof focuses on $\Sigma_+$.  To ease the notations, we write
$dz$ in the place of $dz_+$.  For periodic vortex crystals, we compute the
partial derivatives
\[
  \frac{\partial}{\partial \alpha_+} \int_{A+} \frac{\tilde Q}{dz} = 2\ii,
  \qquad
  \frac{\partial}{\partial \beta_+} \int_{A_+} \frac{\tilde Q}{dz} = 2,
\]
and all other partial derivatives of $\int_{A+} \tilde Q/dz$ vanish.  By
\cite{traizet2008}*{Lemma 3}, we have the partial derivatives
\[
  \frac{\partial \mathcal{E}_k}{\partial t_k} = \frac{\partial}{\partial t_k} \int_{A_k} \frac{(z-p_k)\tilde Q}{dz}
  = \int_{A_k}(z-p_k) \left(  2 \frac{\partial\tphi_1}{\partial t_k}
  -2\ii\frac{\partial\tphi_2}{\partial t_k}\right)
  = -8\pi\ii,
\]
and all other partial derivatives of $\mathcal{E}_k$ vanish.

Finally, we compute, at $\varepsilon = 0$,
\begin{align*}
  (\re \tau, \im \tau) = \re \int_{B_+} (\tphi_1, \tphi_2) = \re \int_0^{\tau_+}(dz, -\ii dz) = (\re \tau_+, \im \tau_+),
\end{align*}
which determines $\tau_+ = \tau$.

The proposition then follows from the Implicit Function Theorem.


\end{document}